\documentclass[a4paper,oneside,11pt,final]{article}

\usepackage[
  a4paper,
  left=20mm,
  right=20mm,
  top=20mm,
  bottom=20mm
]{geometry}

\usepackage[T1]{fontenc}
\usepackage[utf8]{inputenc}

\usepackage{lmodern}

\usepackage{amsmath}
\usepackage{amssymb}
\usepackage{amsthm}
\usepackage{mathtools}

\usepackage{titlesec}

\titleformat{\section}
  {\normalfont\sffamily\Large\bfseries}
  {\thesection}{1em}{}

\titleformat{\subsection}
  {\normalfont\sffamily\large\bfseries}
  {\thesubsection}{1em}{}

\titleformat{\subsubsection}
  {\normalfont\sffamily\normalsize\bfseries}
  {\thesubsubsection}{1em}{}

\titleformat{\paragraph}
  {\normalfont\sffamily\normalsize\bfseries}
  {\theparagraph}{}{}

\usepackage{setspace}
\usepackage[backend=biber,style=numeric,giveninits=true,noerroretextools,maxbibnames=99]{biblatex}
\usepackage{csquotes}

\DeclareSourcemap{
  \maps[datatype=bibtex]{
     \map{
        \step[fieldsource=doi,final]
        \step[fieldset=isbn,null]
        \step[fieldset=issn,null]
        \step[fieldset=url,null]
        }
      \map{
        \step[fieldsource=url,final]
        \step[fieldset=isbn,null]
        \step[fieldset=issn,null]
        }
      \map{
        \step[fieldsource=eprint,final]
        \step[fieldset=url,null]
        }
      }
}

\usepackage{hyperref} 
\usepackage[english]{babel}
\usepackage{enumerate}
\usepackage{aliascnt}
\usepackage{verbatim}
\usepackage[normalem]{ulem}
\usepackage{esint} 

\usepackage[noabbrev,capitalize]{cleveref}

\usepackage[right]{showlabels}

\usepackage[nottoc]{tocbibind} 
\usepackage[title]{appendix}
\usepackage{xcolor}

\usepackage[final]{graphicx}
\usepackage{subcaption}
\usepackage{tikz}
\usetikzlibrary{positioning,arrows.meta,decorations.markings}

\numberwithin{equation}{section} 

\usepackage{hyperref} 
\hypersetup{final,
  colorlinks,
  linkcolor={blue!70!black},
  citecolor={blue!40!black},
  urlcolor={blue!40!black}} 

\let\etoolboxforlistloop\forlistloop 
\usepackage{autonum}
\let\forlistloop\etoolboxforlistloop 

\allowdisplaybreaks[4] 

\newtheoremstyle{sansplain}
  {3pt}   
  {3pt}   
  {\itshape} 
  {}      
  {\sffamily\bfseries} 
  {.}     
  {.5em}  
  {}      

\newtheoremstyle{sansdefinition}
  {3pt}
  {3pt}
  {\itshape}
  {}
  {\sffamily\bfseries}
  {.}
  {.5em}
  {}

\newtheoremstyle{sansremark}
  {3pt}
  {3pt}
  {\normalfont}
  {}
  {\sffamily}
  {.}
  {.5em}
  {}

\theoremstyle{sansplain}

\newtheorem{proposition}{Proposition}[section]

\newaliascnt{lemma}{proposition} 

\newtheorem{lemma}[lemma]{Lemma}

\aliascntresetthe{lemma} 
\Crefname{lemma}{Lemma}{Lemmas}

\newaliascnt{theorem}{proposition} 
\newtheorem{theorem}[theorem]{Theorem}

\aliascntresetthe{theorem} 

\newaliascnt{corollary}{proposition} 

\aliascntresetthe{corollary}

\newaliascnt{hypothesis}{proposition}

\aliascntresetthe{hypothesis} 

\theoremstyle{sansdefinition}

\newaliascnt{definition}{proposition} 
\newtheorem{definition}[definition]{Definition}

\aliascntresetthe{definition} 
\Crefname{definition}{Definition}{Definitions}

\newaliascnt{problem}{proposition} 

\aliascntresetthe{problem} 

\newaliascnt{example}{proposition} 

\aliascntresetthe{example} 

\theoremstyle{sansremark}

\newaliascnt{remark}{proposition} 
\newtheorem{remark}[remark]{Remark}

\aliascntresetthe{remark} 

\def\equationautorefname~#1\null{%
	(#1)\null
}

\newcommand{\R}{\mathbb{R}}

\newcommand{\Z}{\mathbb{Z}}
\newcommand{\N}{\mathbb{N}}
\newcommand{\E}{\mathcal{E}}

\newcommand{\curv}{\vec{\kappa}}
\renewcommand{\S}{\mathbb{S}}

\newcommand{\A}{\mathcal{A}}

\newcommand{\T}{\mathbb{T}}
\renewcommand{\H}{\mathbb{H}}
\newcommand{\W}{\mathcal{W}}

\newcommand*{\dd}{\mathop{}\!\mathrm{d}}

\title{
    Finite time blow-up for the Willmore flow of \linebreak symmetric figure-eight tori
}

\author{
    \large{Fabian Rupp}\thanks{Department of Mathematics, Johannes Kepler University Linz. Altenberger Straße 69, 4040
Linz, Austria.
\texttt{fabian.rupp@jku.at}} 
    \and 
    \large{Manuel Schlierf}\thanks{Department of Mathematics, University of Salzburg, Hellbrunner Str.~34, 5020 Salzburg, Austria. \texttt{math@manuelschlierf.info}}
}

\date{}

\begin{document}
\maketitle
\vspace{-1cm}
\begin{abstract}
    We show that the Willmore flow of symmetric figure-eight tori of revolution produces a finite time singularity, provided the initial energy is below $12\pi$. Our approach relies on the connection with elastic curves in the hyperbolic plane and a fine analysis of Palais--Smale sequences in the relevant energy regime. 
\end{abstract}

{\small
\noindent \textbf{Keywords and phrases:} Willmore flow, rotational symmetry, figure-eight torus, elastic energy, hyperbolic plane, Palais--Smale sequence.

\noindent \textbf{MSC(2020)}: 35B44 (primary); 53E40, 35B40 (secondary)
}

\section{Introduction and main results}

Let $\Sigma$ be a compact surface without boundary and let $f\colon \Sigma\to\R^3$ be an immersion. The \emph{Willmore energy} of $f$ is defined as
\begin{equation}
    \W(f)=\int_{\Sigma} H^2\dd\mu,
\end{equation}
where $\mu$ is the Riemannian measure on $\Sigma$ induced by the pull-back of the Euclidean metric in $\R^3$ and $H = \frac12 (\kappa_1+\kappa_2)$ is the mean curvature. The  functional $\mathcal{W}$ was already studied by Blaschke and Thomsen \cite{blaschkethomsen1929}, who proved its conformal invariance for closed surfaces, but it became more popular by the work of Willmore \cite{willmore1965}, who observed that round spheres are the global minimizers among closed surfaces. He also made the \emph{Willmore conjecture,} stating that the Clifford torus is the minimizer for immersions with $\Sigma = \T^2$. This was proven by Marques--Neves \cite{marquesneves2014}, based on previous works on the existence problem \cite{liyau1982,simon1993,bauerkuwert2003}.

The associated dynamic problem, the \emph{Willmore flow}, is the geometric evolution equation
\begin{align}
\partial_tf = - (\Delta H+|A^0|^2H)\vec N,
\end{align}
for families of immersions $f\colon[0,T)\times\Sigma \to\R^3$, where $\vec N$ is a local unit normal and $A^0$ is the trace-free part of the second fundamental form. This flow was first studied by Simonett \cite{simonett2001} who proved stability of the round sphere. After that, Kuwert--Schätzle  \cite{kuwertschaetzle2002,kuwertschaetzle2001,kuwertschaetzle2004} introduced a curvature-concentration based blow-up procedure that eventually established global existence and convergence for spherical initial data $f_0\colon \S^2\to\R^3$ satisfying $\W(f_0)\leq 8\pi$. This threshold was shown to be sharp by Blatt \cite{blatt2009}, in the following sense: There exist rotational symmetric initial data $f_0\colon \S^2\to\R^3$ with initial energy arbitrarily close to $8\pi$ whose Willmore flow either has finite time blow-up or does not converge as $t\to\infty$. Indeed, rotational symmetry is preserved and would produce a rotational symmetric equilibrium in the limit $t\to\infty$, but there is no critical point that has the correct turning number \cite{bryant1984}. Since this is an indirect argument, it cannot determine if the singularity occurs in finite or infinite time. The same issue arises in the topological obstructions discovered in \cite{maederbaumdicker2025willmoreenergylandscapespheres}.

Blatt's strategy already relies on the rotational symmetry and uses the Bryant--Griffiths formula \cite{bryantgriffiths1986}, relating the  \emph{(hyperbolic) elastic energy} of the profile curve $\gamma$ in the hyperbolic plane $\mathbb{H}^2=\R\times(0,\infty)$ with the Willmore energy of the induced surface of revolution $f_\gamma$ by
\begin{align}\label{eq:BG}
    \W(f_\gamma)=\frac{\pi}{2}\E(\gamma) + 2\pi\chi(\Sigma).
\end{align}
This relation is also valid in the case of closed curves and tori of revolution and was used in \cite{dallacquamullerschatzlespener2020} to show global existence and convergence to the Clifford torus for tori of revolution with initial energy at most $8\pi$. Also here, the threshold is shown to be sharp by a similar argument as in \cite{blatt2009}. In particular, the question about the finiteness of the singular time is left open.

A parametric approach to the Willmore flow was introduced in \cite{palmurellariviere2022,palmurellariviere2024}. However, also in this setting, the question about possible finite time singularities has remained unsolved.
In numerical simulations of Willmore flows, both evidence for finite time blow-up and against it have been observed, see \cite{mayersimonett2002} and \cite{barrettgarckenuernberg2021}, respectively.

In this work, we resolve the question by providing an analytical proof of a finite time singularity.
\begin{theorem}[\Cref{thm:main} below]\label{thm:intro-main}
    Suppose $f_0$ is a torus of revolution such that its profile curve $\gamma_0$ has turning number zero and $\gamma_0(-x) =  \operatorname{diag}(-1,1)\gamma_0(x)$ for $x\in \S^1$.
    If $\W(f_0)\leq 12\pi$, then the Willmore flow starting from $f_{0}$ develops a finite time singularity.
\end{theorem}
The $\lambda$-figure eights that were considered in \cite{dallacquamullerschatzlespener2020} (based on \cite{muellerspener2020}) provide admissible initial data in \Cref{thm:intro-main} and
have energy arbitrarily close to $8\pi$. 

The question about the finiteness of the singular time for the Willmore flow is very challenging for two reasons. First, the Willmore flow is a fourth order parabolic PDE, so no maximum or comparison principles are available. Second, the Willmore energy is scale invariant, making it very difficult to control a priori the behavior of global quantities such as the area or the volume. Except for stationary solutions, we are not aware of any explicit closed solutions, in contrast to other flows like the mean curvature flow where there is a plethora of shrinking or translating solutions. For modifications of the Willmore flow, several works have established the existence of finite time singularities, see, for instance, \cite{mccoywheeler2016,blatt2019,schlierf2025cpde}. However, in all these works, the characteristic scale invariance is broken. This is different for the constrained Willmore flows studied in \cite{jachan2014,rupp2023,rupp2024,ruppscharrerschlierf2024}, where the constraints either determine a fixed scale or are scale-invariant as well.

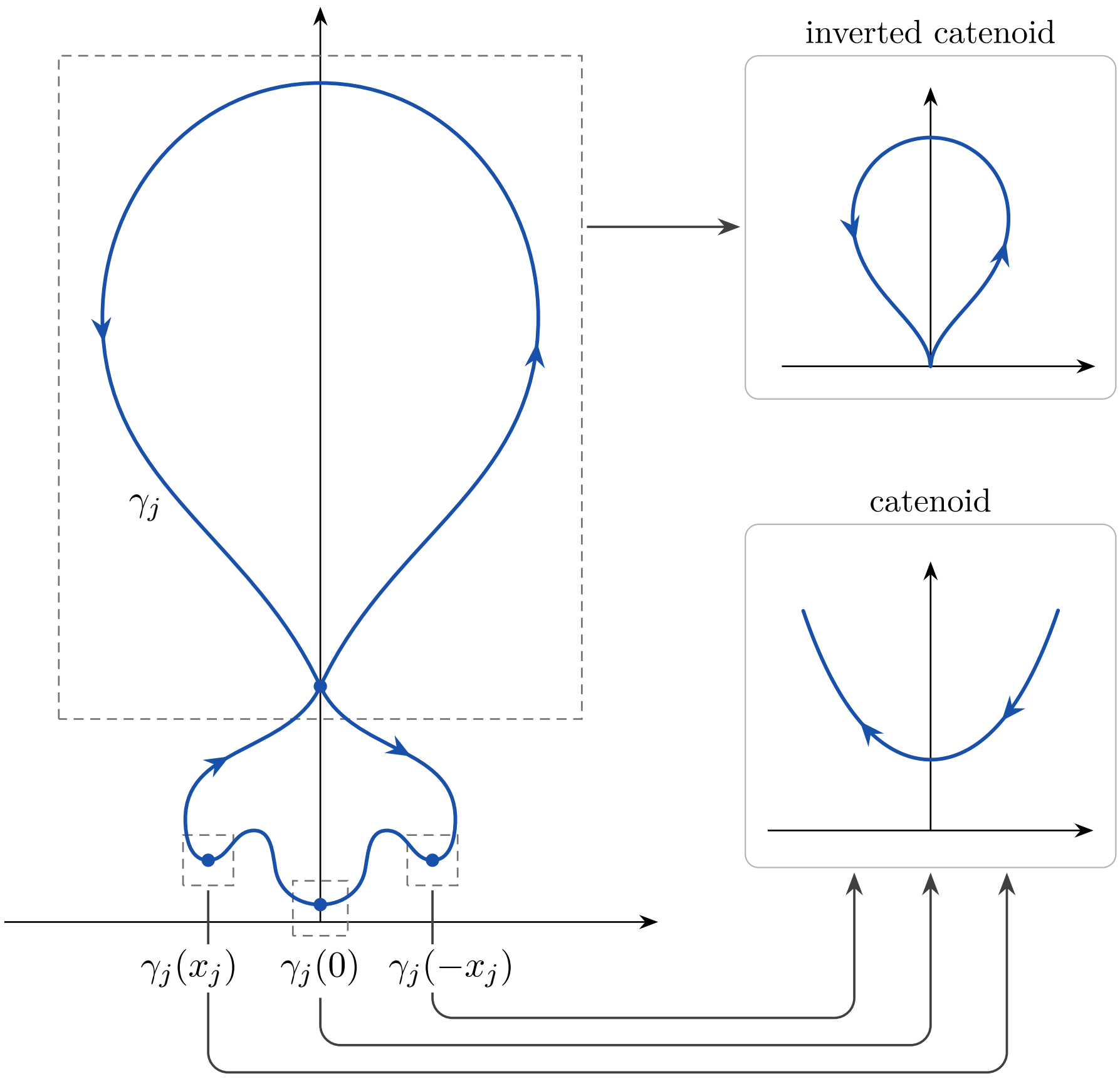
\begin{figure}
\centering
\begin{tikzpicture}[
    scale=1.15,
    profile/.style    = {line width=1.1pt, draw=profilecolor, line cap=round, line join=round},
    oriented/.style 2 args = {postaction={decorate}, decoration={markings,
        mark=at position #1 with {\arrow{Stealth[length=2.6mm,width=2.1mm]}},
        mark=at position #2 with {\arrow{Stealth[length=2.6mm,width=2.1mm]}}}},
    axis/.style       = {line width=0.5pt, -{Stealth[length=2mm]}},
    zoombox/.style    = {densely dashed, line width=0.5pt, draw=black!55},
    inset/.style      = {line width=0.4pt, draw=black!30, rounded corners=4pt},
    blowup/.style     = {line width=0.7pt, draw=black!75, -{Stealth[length=2.4mm]},
                         rounded corners=6pt, shorten <=1.5pt, shorten >=1.5pt},
    insetlabel/.style = {font=\small, above},
    necklabel/.style  = {fill=white, inner sep=1.5pt, anchor=base},
    dot/.style        = {circle, fill=profilecolor, inner sep=0pt, minimum size=4pt},
    declare function  = {
      sech(\t) = 1/cosh(\t);
      invx(\t) = \t*sech(\t)*sech(\t)/(1 + \t*\t*sech(\t)*sech(\t));
      invy(\t) = sech(\t)/(1 + \t*\t*sech(\t)*sech(\t));
    },
  ]
  \definecolor{profilecolor}{RGB}{25,80,170}
  \def\lowscale{0.42}   

  \draw[axis] (-2.9,0) -- (3.1,0);
  \draw[axis] (0,0) -- (0,8.4);
  \coordinate (S) at (0,{5.15*\lowscale});   

  \draw[zoombox] (-2.4,{5.15*\lowscale-0.3}) rectangle (2.4,7.95);
  \draw[profile, oriented={0.24}{0.76}]
    (0,7.7) to[out=180,in=90] (-2.0,5.55) to[out=-90,in=115] (S)
            to[out=65,in=-90] (2.0,5.55)  to[out=90,in=0]    (0,7.7);

  \begin{scope}[scale=\lowscale]
    \draw[zoombox] (-3.0,0.8)  rectangle (-1.9,1.9);    
    \draw[zoombox] (1.9,0.8)   rectangle (3.0,1.9);     
    \draw[zoombox] (-0.6,-0.3) rectangle (0.6,0.9);     
    \draw[profile, oriented={0.14}{0.86}]
      (0,5.15)    to[out=-65,in=90]  (2.95,2.25)  to[out=-90,in=0]   (2.45,1.35)
                  to[out=180,in=0]   (1.45,2.0)   to[out=180,in=80]  (0.98,1.15)
                  to[out=-100,in=0]  (0,0.38)     to[out=180,in=280] (-0.98,1.15)
                  to[out=100,in=0]   (-1.45,2.0)  to[out=180,in=0]   (-2.45,1.35)
                  to[out=180,in=270] (-2.95,2.25) to[out=90,in=245]  (0,5.15);
    \foreach \p in {(-2.45,1.35), (2.45,1.35), (0,0.38)}
      \node[dot] at \p {};
    \coordinate (neckL) at (-2.45,0.8);   
    \coordinate (neckR) at (2.45,0.8);
  \end{scope}
  \node[dot] at (S) {};
  \node[left=3pt] at (-1.24,3.8) {$\gamma_j$};

  \begin{scope}[shift={(5.6,5.1)}, scale=2.1]
    \draw[axis] (-0.65,0) -- (0.72,0);
    \draw[axis] (0,0) -- (0,1.22);
    \draw[profile, oriented={0.25}{0.75}]
      (0,0) -- plot[variable=\t, domain=-4:4, samples=401] ({-invx(\t)}, {invy(\t)}) -- (0,0);
  \end{scope}
  \draw[inset] (3.9,4.8) rectangle (7.3,7.95);
  \node[insetlabel] at (5.6,7.95) {inverted catenoid};
  \draw[blowup] (2.4,6.38) -- (3.9,6.38);

  \begin{scope}[shift={(5.6,0.84)}, scale=0.65]
    \draw[axis] (-2.3,0) -- (2.3,0);
    \draw[axis] (0,0) -- (0,3.8);
    \draw[profile, oriented={0.3}{0.7}, domain=-1.8:1.8, samples=181]
      plot ({-\x}, {cosh(\x)});
  \end{scope}
  \draw[inset] (3.9,0.5) rectangle (7.3,3.65);
  \node[insetlabel] at (5.6,3.65) {catenoid};
  \draw[blowup] (neckR) |- (4.9,-0.88) -- (4.9,0.5);
  \draw[blowup] (neckL) |- (6.3,-1.38) -- (6.3,0.5);

  \node[necklabel, xshift=-2mm] at (neckL |- 0,-0.5) {$\gamma_j(x_j)$};
  \node[necklabel] (x0)         at (0,-0.5)          {$\gamma_j(0)$};
  \node[necklabel, xshift=2mm]  at (neckR |- 0,-0.5) {$\gamma_j(-x_j)$};
  \draw[blowup] (x0.south) |- (5.6,-1.13) -- (5.6,0.5);
\end{tikzpicture}

\caption{A Palais--Smale sequence $\gamma_j$ that cannot occur for $\mathcal{W}\in [8\pi,12\pi)$}\label{fig:introduction}
\end{figure}

We first conjectured that finite time singularities can arise from initial data as in \Cref{thm:intro-main} by the following two observations. On the one hand, if the flow in \Cref{thm:intro-main} existed globally, along some sequence of times, the Willmore gradient would need to converge to zero by the $L^2$-gradient flow structure. On the other hand, in \cite{dallacquamuellerruppschlierf2025}, we showed (under some additional area bounds) that the sequence of profile curves converges to the \emph{inverted catenoid}, see \Cref{fig:introduction}. These two observations seem incompatible, since the inverted catenoid satisfies the weak vectorial Willmore equation, cf.\ \cite{riviere2008}, with a Dirac delta inhomogeneity in direction of the rotation axis \cite{bernardriviere2013}. However, the convergence to the inverted catenoid is only locally smoothly away from the rotation axis, so that this argument by itself is insufficient. This idea seems to be limited to the case of tori and is not directly extendable to the spherical initial data considered in \cite{blatt2009}, see \Cref{rem:spheres_not_good} below.

Nevertheless, our strategy for proving \Cref{thm:intro-main} is heavily inspired by these observations. In \Cref{thm:willgrad-lower-bound}, we prove a scale-invariant bound on the velocity. We show that
\begin{align}\label{eq:intro-velocity_estimate}
    \Vert\nabla \W(f)\Vert_{L^2(\T^2)}^2 \A(f)^2 \geq c>0
\end{align}
for tori of revolution with symmetric figure-eight profile curves in the energy regime determined by \Cref{thm:intro-main}. Here, for an immersion $f\colon \Sigma\to\R^3$, we write $\nabla \W(f) = (\Delta H+|A^0|^2 H)\vec N$ for the $L^2$-gradient of the Willmore energy and $\A(f) = \int_{\Sigma} \dd \mu$ for the total area. An elementary estimate yields that the squared area grows at most affinely in $t$, cf.\ \cite[p.~349]{kuwertschaetzle2004}, yielding
\begin{align}\label{eq:intro_contradiction}
    \Vert \partial_t f\Vert_{L^2(\T^2)}^2 \geq \tilde c (1+t)^{-1}.
\end{align}
However, the energy dissipation gives that $\Vert \partial_t f\Vert_{L^2(\T^2)}^2$ is integrable in $t\in [0,T_{\mathrm{max}})$. Hence $T_{\mathrm{max}}<\infty$ has to hold.

The proof of \eqref{eq:intro-velocity_estimate} (respectively \Cref{thm:willgrad-lower-bound}) uses a contradiction argument that is visualized in \Cref{fig:introduction}. After rescaling, we consider Palais--Smale sequences given by symmetric figure-eight curves $\gamma_j\colon \S^1\to\H^2$. Note that the Palais--Smale condition is not satisfied for the elastica functional in $\mathbb{H}^2$ \cite{linner1998}, see also \cite{langersinger1987,schrader2016}. Indeed, by \cite{muellerspener2020}, the sequence $\gamma_j$ cannot converge smoothly as $j\to\infty$. We give a detailed analysis of its singular limit with the goal of excluding such behavior in the energy regime $\mathcal{W}\in [8\pi, 12\pi)$. We show that $\W(f_{\gamma_j})\to 8\pi$ and that the conformal classes of the immersions $(f_{\gamma_j})_j$ diverge, so this can be seen as a refinement of the energy quantization result in \cite{scharrer2026bubbleclassificationimmersionsboundary} under specific symmetry assumptions.

First, following the arguments in \cite{schlierf2024cvpde,dallacquamuellerruppschlierf2025}, we obtain a macroscopic limit that we identify as the inverted catenoid. Second, and most importantly, we are able to extend the analysis in \cite{dallacquamuellerruppschlierf2025} on the inner part. Around a point on the vertical axis where the \emph{height}, i.e., the distance to the rotation axis, converges to zero, we parametrize by hyperbolic arclength such that $\gamma_j^{(2)}(0)\to 0$, allowing us to find a catenoid profile curve in the blow-up limit, cf.\ \Cref{fig:introduction}. Moreover, we show that the mean curvature of $f_{\gamma_j}$ converges --- without blow-up --- to a solution of the Jacobi equation for the catenoid (viewed as a minimal surface of revolution). The classification of the macroscopic limit combined with an application of Noether's theorem then enables us to establish the existence of two further local minima of the height at $x=\pm x_j$. The Noether quantity we employ arises from the infinitesimal invariance generated by translation along the rotation axis; the same direction where also the first variation of the Willmore energy at the inverted catenoid is nonzero. We show that the three height minima become infinitely far apart (in hyperbolic arclength) and thus the total energy splits. Since each minimum approaches a catenoid, we find that the total energy of the curve must have $\mathcal{W}\geq 12\pi$, closing the contradiction argument.

\section{Preliminaries}

\subsection{Hyperbolic geometry and surfaces of revolution}
Let $\T^2=\S^1\times\S^1$ be the standard $2$-torus with orientation $\{\partial_{x_1},\partial_{x_2}\}$, where $\S^1=\R/(2\pi\Z)$ in the sequel. For $\varphi\colon\S^1\to\R^2$, we define $\mathcal R\varphi\colon\T^2\to\R^3$ by
\begin{equation}\label{eq:def-R}
    (\mathcal R\varphi)(x_1,x_2)\vcentcolon=
    (\varphi^{(1)}(x_1),\varphi^{(2)}(x_1)\cos(x_2),\varphi^{(2)}(x_1)\sin(x_2)).
\end{equation}
For an immersion $\gamma\colon\S^1\to\H^2$, we write $f_\gamma\vcentcolon=\mathcal R\gamma$ for its associated surface of revolution. We use the same notation for curves defined on an interval and for vector fields along a curve. 

On $\H^2=\R\times (0,\infty)$, fix the metric $g_{(y^{(1)},y^{(2)})}\vcentcolon=\frac{1}{(y^{(2)})^2}\langle\cdot,\cdot\rangle$ and denote by $\nabla$ covariant differentiation. 
With $\partial_s\gamma\vcentcolon=\frac{1}{|\partial_x\gamma|_g}\partial_x\gamma$, the curvature of an immersion $\gamma\colon\S^1\to\H^2$ is given by $\curv\vcentcolon=\nabla_{\partial_s\gamma}\partial_s\gamma$. With the hyperbolic unit normal $\vec n_\gamma=(-\partial_s\gamma^{(2)},\partial_s\gamma^{(1)})$, define the \emph{scalar curvature} $\kappa\vcentcolon=\langle\curv,\vec n_\gamma\rangle_g$. For the sake of simplicity, we write $\nabla_s$ for $\nabla_{\partial_s\gamma}$, and for a vector field $V$ along $\gamma$, let $V^\perp$ denote its projection onto the subspace orthogonal to $\partial_s\gamma$, and set $\nabla_s^\perp V\vcentcolon=(\nabla_sV)^\perp$. In particular,
\begin{equation}\label{eq:frenet}
    \nabla_s\partial_s\gamma=\kappa\vec n_\gamma, \quad \nabla_s\vec n_\gamma=-\kappa\partial_s\gamma, \quad\text{and}\quad \nabla_s^\perp\vec n_\gamma=0.
\end{equation}
Moreover, cf.\ \cite[Equation (2.5)]{dallacquamuellerruppschlierf2025}, we have for a vector field $V\colon\S^1\to\R^2$ along $\gamma$ that
\begin{equation}\label{eq:cov-der-H2}
    \nabla_s V = \partial_s V - \frac{V^{(2)}}{\gamma^{(2)}}\partial_s\gamma + \frac{V^{(1)}}{\gamma^{(2)}}\vec n_\gamma.
\end{equation}
Writing $\dd s\vcentcolon=|\partial_x\gamma|_g\dd x$ for the arclength element, the \emph{hyperbolic length} and the \emph{elastic energy} are
\begin{equation}
    \mathcal L(\gamma)\vcentcolon=\int_{\S^1}\dd s
    \quad\text{and}\quad
    \E(\gamma)\vcentcolon=\int_{\S^1}|\curv|_g^2\dd s.
\end{equation}
We also use the Euclidean length $\mathcal L_{\R^2}(\gamma)\vcentcolon=\int_{\S^1}|\partial_x\gamma|\dd x$, arclength derivative $\partial_{s}^{\R^2}=\frac1{|\partial_x\gamma|}\partial_x=\frac1{\gamma^{(2)}}\partial_s$, and the turning number 
\begin{equation}
    \mathcal T(\gamma)=\frac{1}{2\pi}\int_{\S^1}\kappa_{\R^2}|\partial_x\gamma|\dd x\in\Z,
\end{equation}
where $\kappa_{\R^2} = \langle (\partial_{s}^{\R^2})^2\gamma,\vec n_{\gamma}^{\R^2}\rangle$ is the Euclidean scalar curvature and $\vec n_{\gamma}^{\R^2} = (-\partial_{s}^{\R^2}\!\gamma^{(2)},\partial_{s}^{\R^2}\!\gamma^{(1)})$.

For the surface of revolution $f_\gamma$, the induced metric and unit normal are 
\begin{equation}\label{eq:surf-rev-metric-and-normal}
    g_{f_\gamma}=|\partial_x\gamma(x_1)|^2\dd x_1^2+(\gamma^{(2)}(x_1))^2\dd x_2^2\quad\text{and}\quad\vec N_{f_\gamma}=-(\gamma^{(2)})^{-1}\mathcal R\vec n_\gamma=-\mathcal R\vec n_\gamma^{\R^2}\!.
\end{equation}
In particular, 
\begin{equation}\label{eq:area-measure}
    \dd\mu_{f_\gamma} = \gamma^{(2)}(x_1)|\partial_x\gamma(x_1)|\dd x_1\otimes \dd x_2 = (\gamma^{(2)}(x_1))^2\dd s \otimes \dd x_2.
\end{equation}

\begin{proposition}\label{prop:elastic-willmore-identities}
    Let $I=[a,b]$ or $I=\S^1$, and let $\gamma\colon I\to\H^2$ be a smooth immersion. For $I=\S^1$, the boundary terms below are omitted. For the $L^2(\dd s)$-gradient of $\E$ given by
    \begin{equation}\label{eq:geometry}
        \nabla\E(\gamma)=2(\nabla_s^\perp)^2\curv+|\curv|_g^2\curv-2\curv
    \end{equation}
    and for every smooth vector field $V\colon I\to\R^2$ along $\gamma$, we have
    \begin{align}
        \left.\frac{\dd}{\dd t}\E(\gamma+tV)\right|_{t=0} &= \int_I\langle\nabla\E(\gamma),V\rangle_g\dd s  + 2\Big[\langle\curv,\nabla_sV\rangle_g-\langle\nabla_s^\perp\curv,V\rangle_g-\frac12|\curv|_g^2\langle V,\partial_s\gamma\rangle_g\Big]_a^b.\label{eq:elastic-first-variation}
    \end{align}
    Moreover, for the surface of revolution $f_\gamma=\mathcal R\gamma$, the $L^2$-gradients of $\E$ and $\W$ are related by
    \begin{equation}\label{eq:gradient}
        \nabla\W(f_\gamma)=\mathcal R\Big(\frac{1}{4(\gamma^{(2)})^4}\nabla\E(\gamma)\Big)
    \end{equation}
    so that, in particular,
    \begin{equation}\label{eq:gradient-norm}
        \|\nabla\W(f_\gamma)\|_{L^2(\dd\mu_{f_\gamma})}^2=\frac{\pi}{8}\int_I\frac{|\nabla\E(\gamma)|_g^2}{(\gamma^{(2)})^4}\dd s.
    \end{equation}

    We call $\gamma$ an \emph{elastica} if it is critical for $\E$ with respect to compactly supported variations away from $\partial I$. This is equivalent to $f_\gamma$ being a \emph{Willmore immersion}.
\end{proposition}
\begin{proof}
    The computation in \cite[Lemma~2.4 and Remark~2.5]{dallacquaspener2017}, keeping track of the boundary terms, yields \eqref{eq:elastic-first-variation}. For the identity \eqref{eq:gradient}, see \cite[Theorem~4.1]{dallacquaspener2018} and \eqref{eq:surf-rev-metric-and-normal}.

    For a vector field $X$ along $\gamma$, we have $|\mathcal RX|_{\R^3}^2=(\gamma^{(2)})^2|X|_g^2$. Combining \eqref{eq:area-measure} with \eqref{eq:gradient} yields \eqref{eq:gradient-norm}. The characterization of elastica follows from \eqref{eq:elastic-first-variation}, \eqref{eq:geometry}, and \eqref{eq:gradient}.
\end{proof}

\begin{definition}[Noether quantity]\label{rem:elastic-noether}
    Let $V$ be the restriction to $\gamma$ of a Killing field on $(\H^2,g)$. We define
    \begin{equation}\label{eq:elastic-noether-quantity}
        P_{\gamma,V}\vcentcolon=2\langle\nabla_s^\perp\curv,V\rangle_g+|\curv|_g^2\langle V,\partial_s\gamma\rangle_g-2\langle\curv,\nabla_sV\rangle_g \colon I\to\R.
    \end{equation}
    In particular, for $V=(1,0)$, by \eqref{eq:frenet} and \eqref{eq:cov-der-H2} we obtain
    \begin{align}\label{eq:horizontal-noether}
        P_\gamma &\vcentcolon= P_{\gamma,(1,0)} = \frac{1}{\gamma^{(2)}} \Big(-2\partial_s\kappa\, \partial_s^{\R^2}\!\gamma^{(2)}+\kappa^2\partial_s^{\R^2}\!\gamma^{(1)}
        -2\kappa\Big).
    \end{align}
\end{definition}

Note that $-P_{\gamma,V}$ is the boundary integrand in \eqref{eq:elastic-first-variation}. Since $V$ is a Killing field, we have $\partial_{t,0}\E(\gamma+tV|_J)=0$ for any sub-interval $J\subset I$, and thus \eqref{eq:elastic-first-variation} implies $
    \partial_sP_{\gamma,V}=\langle\nabla\E(\gamma),V\rangle_g$.
In particular, $P_{\gamma,V}$ is constant along every elastica. This can be viewed as the \emph{conservation law} associated with $V$ by Noether's theorem, cf.\ \cite[Theorem~(2.a.8)]{bryantgriffiths1986}, see also \cite{langersinger1984,bernard2016}. The Killing field $V=(1,0)$ is the infinitesimal generator of horizontal translations. This particular choice has also been used for the elastica functional in \cite{dallacquadeckelnickwheeler2013}. 
As a quantitative version of Noether's theorem for elastica we have the estimate
\begin{align}
        |\partial_sP_\gamma|=|\langle\nabla\E(\gamma),(1,0)\rangle_g|\leq \frac{|\nabla\E(\gamma)|_g}{\gamma^{(2)}}.\label{eq:horizontal-noether_estimate}
\end{align}

Next, we relate the curvatures of the profile curves and the associated surfaces of revolution.
\begin{proposition}\label{rem:mean-curvature-identities}
    Let $\gamma\colon I\to\H^2$ be a smooth immersion. Then, we have
    \begin{equation}\label{eq:H-A0}
        H_{f_\gamma}=\frac{2\partial_s^{\R^2}\!\gamma^{(1)}-\kappa}{2\gamma^{(2)}}\quad\text{and}\quad |A_{f_\gamma}^0|^2=\frac{\kappa^2}{2(\gamma^{(2)})^2},
    \end{equation}
    and $H_{f_\gamma}$ is a solution to the inhomogeneous ODE
    \begin{align}
        2\partial_s^2H_{f_\gamma}+\kappa^2H_{f_\gamma} 
        &= -\frac{1}{2\gamma^{(2)}}\langle\nabla\E(\gamma),\vec n_\gamma\rangle_g.\label{eq:forced-Jacobi-equation-H}
    \end{align}
    Moreover, the Noether quantity $P_\gamma$ defined in \eqref{eq:horizontal-noether} satisfies 
    \begin{equation}\label{eq:horizontal-noether-Hformula}
        \frac{P_\gamma}{4}=\partial_s^{\R^2}\!\gamma^{(2)}\partial_sH_{f_\gamma}+\gamma^{(2)}\,\partial_s^{\R^2}\!\gamma^{(1)}H_{f_\gamma}^2-(\partial_s^{\R^2}\!\gamma^{(1)})^2H_{f_\gamma}.
    \end{equation}
\end{proposition}
\begin{proof}
    For \eqref{eq:H-A0}, using \eqref{eq:surf-rev-metric-and-normal} and \cite[Example~4 on p.~161]{docarmo1976} (with opposite orientation) we have that the principal curvatures of $f_\gamma$ are given by
    \begin{equation}
        \kappa_{f_\gamma,1}=-\kappa_{\R^2}\quad\text{and}\quad  \kappa_{f_\gamma,2}=\frac{\partial_s^{\R^2}\!\gamma^{(1)}}{\gamma^{(2)}}.
    \end{equation}
    Using \cite[Equation~(2.2)]{eichmanngrunau2019}, we obtain 
    \begin{equation}
        \kappa = \gamma^{(2)} \kappa_{\R^2} + \partial_s^{\R^2}\!\gamma^{(1)}.
    \end{equation}
    Altogether, since $H_{f_\gamma}=\frac12(\kappa_{f_\gamma,1}+\kappa_{f_\gamma,2})$ and $|A^0_{f_\gamma}|^2=\frac12(\kappa_{f_\gamma,1}-\kappa_{f_\gamma,2})^2$, \eqref{eq:H-A0} follows. Using that $H_{f_\gamma}$ does not depend on $x_2$, this yields
    \begin{align}
        2\partial_s^2H_{f_\gamma}+\kappa^2H_{f_\gamma} &= 2(\gamma^{(2)})^2\big(\Delta_{g_{f_\gamma}}H_{f_\gamma}+|A_{f_\gamma}^0|^2H_{f_\gamma}\big) = 2(\gamma^{(2)})^2\langle\nabla\W(f_\gamma),\vec N_{f_\gamma}\rangle
    \end{align}    
    which proves \eqref{eq:forced-Jacobi-equation-H} when combined with \eqref{eq:gradient} and \eqref{eq:surf-rev-metric-and-normal}.
    Finally, differentiating the equation for $H_{f_\gamma}$ and using \eqref{eq:frenet} and \eqref{eq:cov-der-H2}, we obtain
    \begin{equation}\label{eq:ds-H}
        \partial_sH_{f_\gamma}=-\frac{\partial_s\kappa+\kappa\, \partial_s^{\R^2}\!\gamma^{(2)}}{2\gamma^{(2)}}. 
    \end{equation}
    Together with \eqref{eq:horizontal-noether} and \eqref{eq:H-A0}, we get \eqref{eq:horizontal-noether-Hformula}.
\end{proof}

\subsection{Hyperbolic elastica}
Using \eqref{eq:frenet} and \eqref{eq:geometry}, an immersion $\gamma$ is an elastica if and only if its hyperbolic curvature satisfies
\begin{equation}\label{eq:free-elastica-ode}
    2\partial_s^2\kappa+\kappa^3-2\kappa=0.
\end{equation}
We say that $\gamma$ is parametrized by hyperbolic arclength if $|\partial_x\gamma|_g\equiv 1$. In this case, \eqref{eq:free-elastica-ode} turns into an ODE for the function $s\mapsto \kappa(s)$, and all possible solutions are classified in \cite{langersinger1984}, also see \cite[Proposition~2.8]{muellerspener2020}. We use the following consequence of this classification.

\begin{proposition}\label{lem:finite-energy-elastica}
    Let $I\subset\R$ be an interval of infinite length, and let $I\ni s\mapsto \kappa(s)\in\R$ be a $W_{\mathrm{loc}}^{2,2}$-solution of the ODE \eqref{eq:free-elastica-ode}. If $\int_{I}\kappa^2(s)\dd s<\infty$, then either $\kappa\equiv0$, or there exist $a\in\R$, $\sigma\in\{-1,1\}$ such that $\kappa(s)= 2\sigma/\cosh(a+s)$ for all $s\in I$ and $\int_{I}\kappa^2(s)\dd s=4[\tanh(a+s)]_{\partial I}$. 
\end{proposition}
Curves $\gamma$ corresponding to the second case in \Cref{lem:finite-energy-elastica} are called \emph{asymptotically geodesic elastica.}
\begin{proof}[Proof of \Cref{lem:finite-energy-elastica}]
    Smoothness of $\kappa$ is immediate from \eqref{eq:free-elastica-ode}. By \cite[Proposition~2.8]{muellerspener2020}, every smooth solution of \eqref{eq:free-elastica-ode} other than zero and translations of $\pm 2/\cosh(s)$ has nonzero periodic curvature, and thus has infinite $L^2$-norm on $I$. If $\kappa$ is a translation of $\pm 2/\cosh(s)$, the choice of the sign $\pm$ is constant, and we may use $\partial_s\tanh = 1/\cosh^2$ to compute the $L^2$-norm.
\end{proof}
We next recall a condition for the associated curves distinguishing the two cases $\kappa\equiv 0$ and $\kappa(s)= 2\sigma/\cosh(a+s)$ in \Cref{lem:finite-energy-elastica}, see \cite[Remark~2.17]{schlierf2024cvpde} and \cite[Remark~5.1]{dallacquamuellerruppschlierf2025}.
If $I=\S^1$ or $I\subset\R$ is an open interval, we call an immersion $\gamma\colon I\to\H^2$ \emph{complete} if $(I,\gamma^*g_{\H^2})$ is complete. 

\begin{proposition}\label{rem:complete-elastica-ends}
If $\gamma$ is a complete elastica with $\E(\gamma)<\infty$ whose two ends converge to the same point in $\R\times\{0\}$, then $\gamma$ is asymptotically geodesic, $\E(\gamma)=8$,  and $f_\gamma$ is an inverted catenoid. Up to isometries of $\H^2$ and reparametrization, $\gamma$ is given by  $s\longmapsto\frac{(s,\cosh(s))}{s^2+\cosh^2(s)}$. In this case, the mean curvature $H$ of the associated surface of revolution is given by 
    \begin{equation}\label{eq:H-inverted-catenoid}
        H(s)=2(1-s\tanh (s))\quad\text{for all $s\in\R$}.
    \end{equation}
\end{proposition}
\begin{proof}
    This follows directly from \Cref{lem:finite-energy-elastica}, \cite[Equation (A.6)]{schlierf2024cvpde}, and \eqref{eq:H-A0}.
\end{proof}

Using one-dimensional elliptic and interpolation estimates, the ODE rigidity in \Cref{lem:finite-energy-elastica} also characterizes sequences $(\kappa_j)_j$ for which \eqref{eq:free-elastica-ode} is only satisfied in the limit.

\begin{lemma}\label{prop:kappaj-palais-smale-rigidity}
    Let $\kappa_j\colon\R\to\R$ be $L_j$-periodic $W^{2,2}_{\mathrm{loc}}$-functions, let $I_j=[-L_j/2,L_j/2)$ and suppose that $L_j\to\infty$. Assume for $j\to\infty$
    \begin{equation}\label{eq:kappaj-palais-smale-assumptions}
        \|\kappa_j\|^2_{L^2(I_j)} \to E\in (0,\infty)\quad\text{and}\quad \|2\partial_s^2\kappa_j+\kappa_j^3-2\kappa_j\|_{L^2(I_j)}\to 0.
    \end{equation}
    Then   
    \begin{equation}\label{eq:kappaj-H2-bounds}
        \sup_{j\in\N}\big(\|\kappa_j\|_\infty+\|\partial_s\kappa_j\|_{L^2(I_j)} + \|\partial_s^2\kappa_j\|_{L^2(I_j)}\big)<\infty,
    \end{equation}
    and there exists $m\in\N$ such that $E=8m$. Let $\operatorname{dist}_j$ denote the distance on $\R/(L_j\Z)$. For all $1\leq i\leq m$ there are $a_j^{(i)} \in \R$, $j\in\N$,  with $\operatorname{dist}_j(a_j^{(i_1)},a_j^{(i_2)})\to\infty$ for $i_1\neq i_2$, and there are $\sigma_i\in\{\pm 1\}$ such that, up to subsequence,
    \begin{equation}\label{eq:kappaj-convergence-at-a_ji}
        \kappa_j(a_j^{(i)}+\cdot)\to \frac{2\sigma_i}{\cosh(\cdot)}\quad\text{in $C^{1,\beta}_{\mathrm{loc}}(\R)$ for all $0<\beta<1/2$ and for all $1\leq i\leq m$}.
    \end{equation}
    In particular, for any $J_j\subset \R$ such that $\operatorname{dist}_j(J_j,a_j^{(i)})\to\infty$ for all $1\leq i\leq m$, we have 
    \begin{equation}\label{eq:kappaj-univ-convergence-to-zero}
        \|\kappa_j\|_{L^\infty}(J_j) \to 0.
    \end{equation}
    Moreover, the sequences $(a_j^{(i)})_{j\in\N}$, $1\leq i\leq m$, are uniquely determined in the following sense: If $(a_j)_{j}$ is a sequence with $\kappa_j(a_j)\to \pm 2$, then 
    \begin{equation}\label{eq:aij-uniqueness}
        \min_{1\leq i\leq m} \mathrm{dist}_j(a_j,a_j^{(i)})\to 0\quad\text{for $j\to\infty$}.
    \end{equation}
\end{lemma}
This result can be seen as a stationary version of \cite[Theorem~1.1]{schlierf2024cvpde}, replacing the parabolic estimates along the hyperbolic elastic flow with the Palais--Smale condition \eqref{eq:kappaj-palais-smale-assumptions}. In fact, the rigidity $E=8m$ we prove in \Cref{prop:kappaj-palais-smale-rigidity} improves \cite[Equation~(1.2)]{schlierf2024cvpde}, which reads as $E\geq 8m$ in our notation.

\begin{proof}[Proof of \Cref{prop:kappaj-palais-smale-rigidity}]
    By Gagliardo--Nirenberg on $I_j$ (cf.\ \cite[Lemma 2.4]{dziukkuwertschaetzle2002}), we get
    \begin{equation}\label{eq:scalar-interpolation}
        \Vert\kappa_j^3\Vert_{L^2(I_j)}=\Vert\kappa_j\Vert_{L^6(I_j)}^3\leq C(\Vert\partial_s^2\kappa_j\Vert_{L^2(I_j)}^{1/2}\Vert\kappa_j\Vert_{L^2(I_j)}^{5/2}+1/L_j\|\kappa_j\|_{L^2(I_j)}^3) \leq C(1+\Vert\partial_s^2\kappa_j\Vert_{L^2(I_j)}^{1/2}\Vert\kappa_j\Vert_{L^2(I_j)}^{5/2}).
    \end{equation}
    Here and in the following, $C>0$ denotes a constant which may be taken independent of $j$ using that $|I_j|=L_j\to\infty$, so $\sup_j 1/L_j$ is bounded. Using \eqref{eq:kappaj-palais-smale-assumptions}, Young's inequality yields
    \begin{align}
        2\Vert\partial_s^2\kappa_j\Vert_{L^2(I_j)}&\leq\Vert2\partial_s^2\kappa_j+\kappa_j^3-2\kappa_j\Vert_{L^2(I_j)}+2\Vert\kappa_j\Vert_{L^2(I_j)}+\Vert\kappa_j^3\Vert_{L^2(I_j)} \leq C+C\Vert\partial_s^2\kappa_j\Vert_{L^2(I_j)}^{1/2}\\
        &\leq C+\Vert\partial_s^2\kappa_j\Vert_{L^2(I_j)}.\label{eq:scalar-h2-bound}
    \end{align}
    Again, all constants are independent of $j$. 
    Periodicity and integration by parts imply \eqref{eq:kappaj-H2-bounds}. Thus, for any sequence $(a_j)_j$ of real numbers, $(\kappa_j(a_j+\cdot))_j$ is bounded in $W^{2,2}_{\mathrm{loc}}(\R)$. If $\liminf_{j\to\infty}|\kappa_j(a_j)|>0$, using \Cref{lem:finite-energy-elastica} and \eqref{eq:kappaj-palais-smale-assumptions}, after passing to a subsequence, there exist $a_*\in\R$ and $\sigma\in\{-1,1\}$ such that
    \begin{align}
        \kappa_j(a_j+\cdot) &\to \frac{2\sigma}{\cosh(a_*+\cdot)}\text{ in $C^{1,\beta}_{\mathrm{loc}}(\R)$ for all $0<\beta<\frac12$, and}\\
        \quad \lim_{j\to\infty} & \|\kappa_j\|_{L^2(a_j-R,a_j+R)}^2 = 8+o(1)\text{ as $R\to\infty$}.\label{eq:compactness-for-selected-aj}
    \end{align}
    Let $J_j\subset I_j$ be intervals with length at least 1. Testing the left-hand side in \eqref{eq:free-elastica-ode} (for $\kappa=\kappa_j$) with $\kappa_j$, integrating by parts, estimating the boundary terms with $\|\partial_s\kappa_j\|_{L^{\infty}(J_j)}\leq C\|\kappa_j\|_{W^{2,2}(J_j)}^{\frac34}\|\kappa_j\|_{L^2(J_j)}^{\frac14}\leq C$ by Gagliardo--Nirenberg and the above estimates, and using \eqref{eq:kappaj-palais-smale-assumptions}, we get
    \begin{align}
        2\Vert\partial_s\kappa_j\Vert_{L^2(J_j)}^2+2\Vert\kappa_j\Vert_{L^2(J_j)}^2&=\Vert\kappa_j\Vert_{L^4(J_j)}^4 + 2[\partial_s\kappa_j\kappa_j]_{\partial J_j} - \int_{J_j}(2\partial_s^2\kappa_j+\kappa_j^3-2\kappa_j)\kappa_j\dd s
        \\
        &\leq \|\kappa_j\|_{L^\infty(J_j)}(\|\kappa_j\|_{L^\infty(J_j)}\|\kappa_j\|_{L^2(J_j)}^2 + C) + o(1)\\
        &\leq C \|\kappa_j\|_{L^\infty(J_j)} + o(1)\text{ as $j\to\infty$},\label{eq:scalar-testing}
    \end{align}
    using $\|\kappa_j\|_{L^\infty(J_j)}\leq C$ by \eqref{eq:kappaj-H2-bounds} in the last estimate. Consequently, $\|\kappa_j\|_{L^\infty(J_j)}\to 0$ implies that $\Vert\kappa_j\Vert_{L^2(J_j)}\to 0$.

    Setting $J_j=I_j$, $E>0$ thus yields that there exist $h_j\in I_j$ with $|\kappa_j(h_j)|\not\to 0$. So applying \eqref{eq:compactness-for-selected-aj}, setting $a_j^{(1)}\in I_j$ as the representative of $h_j-a_*$ mod $L_j$, we get $E\geq 8$ and \eqref{eq:kappaj-convergence-at-a_ji} for $i=1$.
    
    We argue by induction. Suppose that we have constructed sequences $a_j^{(i)}\in I_j$ for $1\leq i\leq k$ satisfying \eqref{eq:kappaj-convergence-at-a_ji} and $\operatorname{dist}_j(a_j^{(i_1)},a_j^{(i_2)})\to\infty$ for $i_1\neq i_2$. By \eqref{eq:compactness-for-selected-aj}, we have $E\geq8k$. If $E>8k$, choose $R_j\to\infty$ sufficiently slowly such that the intervals $B_j^{(i)}\vcentcolon=\{s\in \R/(L_j\Z) \mid \operatornamewithlimits{dist}_j(s,a_j^{(i)})<R_j\}$ are pairwise disjoint, the connected components $J_j^{(1)},\dots,J_j^{(k)}$ of $\R/(L_j\Z)\setminus \bigcup_{i=1}^k B_j^{(i)}$ all have length at least 1, and 
    \begin{equation}
        \sum_{i=1}^k\|\kappa_j\|_{L^2(B_j^{(i)})}^2\to8k.
    \end{equation}
    By \eqref{eq:scalar-testing}, we obtain as $j\to\infty$ that
    \begin{equation}
        E-8k+o(1)=\sum_{i=1}^k\|\kappa_j\|_{L^2(J_j^{(i)})}^2\leq C\sum_{i=1}^k\|\kappa_j\|_{L^\infty(J_j^{(i)})}+o(1).
    \end{equation}
    Thus, up to subsequence, there exist $h_j\notin\bigcup_{i=1}^k B_j^{(i)}$ satisfying $\liminf_{j\to\infty}|\kappa_j(h_j)|>0$ as well as $\operatorname{dist}_j(h_j,a_j^{(i)})\geq R_j$ for all $1\leq i\leq k$. Applying \eqref{eq:compactness-for-selected-aj} and setting $a_j^{(k+1)}=h_j-a_*$ modulo $L_j$, we obtain \eqref{eq:kappaj-convergence-at-a_ji} and $\operatorname{dist}_j(a_j^{(k+1)},a_j^{(i)})\to\infty$ for all $1\leq i\leq k$. In particular, $E\geq8(k+1)$. Since $E<\infty$, the induction terminates with $E=8m$ for some $m\in\N$.

    Finally, if there were $a_j\in\R$ with $\operatorname{dist}_j(a_j,a_j^{(i)})\to\infty$ for all $1\leq i\leq m$ and $|\kappa_j(a_j)|\not\to0$, then, up to subsequence, \eqref{eq:compactness-for-selected-aj} would imply $E\geq8(m+1)$, a contradiction. This proves \eqref{eq:kappaj-univ-convergence-to-zero}. In particular, this applies to any choice of $a_j\in I_j$, $j\in\N$, with $|\kappa_j(a_j)|\to 2$. Combined with $\operatorname{dist}_j(a_j^{(i_1)},a_j^{(i_2)})\to\infty$ for $i_1\neq i_2$, there exist $1\leq i\leq m$ and $(N_j)_j\subset\Z$ such that, up to subsequence, $a_j+N_jL_j-a_j^{(i)}\to a_*\in\R$. Then \eqref{eq:kappaj-convergence-at-a_ji} yields
    \begin{equation}
        \kappa_j(a_j)=\kappa_j(a_j^{(i)}+(a_j+N_jL_j-a_j^{(i)})) \to \frac{2\sigma_i}{\cosh(a_*)}.
    \end{equation}
    So $\cosh(a_*)=1$ which implies $a_*=0$ and thus \eqref{eq:aij-uniqueness}.
\end{proof}

\subsection{Compactness for symmetric Palais--Smale sequences with zero turning number}

Throughout this subsection, we consider a sequence $\gamma_j\colon\S^1\to\H^2$ of smooth immersions satisfying 
\begin{equation}\label{eq:symmetry-turningnumber}
    \gamma_j(-x)=\big(\begin{smallmatrix}-1&0\\0&1\end{smallmatrix}\big)\,\gamma_j(x)\text{ for all $x\in\S^1$} \quad\text{and}\quad  \mathcal T(\gamma_j)=0
\end{equation}
for all $j\in\N$, and suppose that there exists $\delta\in(0,8)$ with
\begin{equation}\label{eq:compactness-assumptions}
    \E(\gamma_j)\leq M=24-\delta,\quad \A(f_{\gamma_j})=1\text{ for all $j\in\N$}, \quad\text{and}\quad \|\nabla\W(f_{\gamma_j})\|_{L^2(\dd\mu_{f_{\gamma_j}})}\to0.
\end{equation}
This subsection is split into three lemmas, each concerned with a different aspect of the asymptotic behavior of this fixed sequence $(\gamma_j)_j$.

The first of the three results applies the compactness and singular-limit arguments in \cite[Lemma~3.4 and proof of Theorem~1.1]{schlierf2024cvpde} and \cite[Proposition~5.4]{dallacquamuellerruppschlierf2025} to the sequence $(\gamma_j)_j$, replacing the estimates along the Willmore flow by \eqref{eq:compactness-assumptions}. While the precise statement looks quite technical, one can read it as follows: Away from the rotation axis, one has local convergence in $C^{3,\beta}$ to an inverted catenoid. This allows computing limits of the hyperbolic curvature $\kappa$ and the mean curvature $H$, after choosing suitable reparametrizations.

\begin{lemma}[Macroscopic convergence to the inverted catenoid]\label{lem:compactness-inverted-catenoid}
    We have $\mathcal L(\gamma_j)\to\infty$ and, up to passing to a subsequence and if necessary reparametrizing each $\gamma_j$ by $x\mapsto \pm x+\pi$ or $x\mapsto \pm x$, we have $\partial_s^{\R^2}\!\gamma_j(\pi)=(-1,0)$ for all $j\in\N$ and
    \begin{align}
        \gamma_j^{(2)}(0)\to 0 &\quad\text{and there exists $\rho_\infty \in (0,\infty)$ with}\quad \gamma_j^{(2)}(\pi)\to \rho_\infty.\label{eq:height-limits-at-0-pi}
    \end{align}
    Denote by $\tilde\gamma_{j,\pi}\colon\R\to\H^2$ the $\mathcal L(\gamma_j)$-periodic lift of the hyperbolic arclength reparametrization of $\gamma_j$ with $\tilde\gamma_{j,\pi}(0)=\gamma_j(\pi)$, for all $j\in\N$. Each $\tilde\gamma_{j,\pi}$ satisfies \eqref{eq:symmetry-turningnumber}, and after passing to a further subsequence, we have that $\tilde\gamma_{j,\pi}\to \tilde\gamma_{\infty,\pi}$ in $C^{3,\beta}_{\mathrm{loc}}(\R)$ for all $0<\beta<1/2$ where
    \begin{equation}\label{eq:gamma-infty-pi}
        \tilde\gamma_{\infty,\pi}(s)= \rho_\infty \frac{(-s,\cosh(s))}{s^2+\cosh^2(s)}\quad\text{for all $s\in\R$},
    \end{equation}
    i.e., $f_{\tilde\gamma_{\infty,\pi}}$ is an inverted catenoid.
    In particular, for the curvatures $\tilde\kappa_{j,\pi}$ of $\tilde\gamma_{j,\pi}$, we have
    \begin{equation}\label{eq:kappa-infty-pi}
        \tilde\kappa_{j,\pi}\to \frac{2}{\cosh(\cdot)}\quad \text{in $C^{1,\beta}_{\mathrm{loc}}(\R)$ for every $0<\beta<1/2$, and thus}\quad P_{\gamma_j}(\pi)\to -\frac{8}{\rho_\infty}.
    \end{equation}
    Furthermore, denoting by $\tilde H_{j,\pi}$ the mean curvature of $f_{\tilde \gamma_{j,\pi}}$, for every $0<\beta<1/2$ we have
    \begin{equation}\label{eq:partc-Hinfty}
    \tilde H_{j,\pi}\to\tilde H_{\infty} \text{ in $C^{1,\beta}_{\mathrm{loc}}(\R)$, where }\quad
        \tilde H_{\infty}(s)=-\frac{2}{\rho_\infty} (1-s \tanh(s))\text{ for all $s\in\R$}.
    \end{equation}
\end{lemma}
\begin{proof}
    By \cite[Lemma~5.3 and proof of Proposition~5.4, Claim~1]{dallacquamuellerruppschlierf2025}, $\A(f_{\gamma_j})=1$ implies
    \begin{equation}\label{eq:euclidean-bounds}
        0<c\leq\max_{\S^1}\gamma_j^{(2)}\leq C\quad\text{and}\quad 0<c\leq\mathcal L_{\R^2}(\gamma_j)\leq C,
    \end{equation}
    with $c,C>0$ independent of $j$. Combining \eqref{eq:euclidean-bounds} with \eqref{eq:gradient-norm}, we obtain
    \begin{equation}\label{eq:elastic-gradient-small}
        \int_{\S^1}|\nabla\E(\gamma_j)|_g^2\dd s_j\leq C\|\nabla\W(f_{\gamma_j})\|_{L^2(\dd\mu_{f_{\gamma_j}})}^2\to0.
    \end{equation}
    In particular, we note that reparametrizing by hyperbolic arclength and using \eqref{eq:frenet} and $\mathcal L(\gamma_j)\geq 4\pi^2/M$ by \cite[Equation~(3.5)]{dallacquamuellerruppschlierf2025}, the same argument as used for \eqref{eq:kappaj-H2-bounds} yields
    \begin{equation}\label{eq:curvature-bounds}
        \sup_j\Big( \|\curv_j\|_{L^\infty} + \|\nabla_s^\perp\curv_j\|_{L^2(\dd s_j)} + \|(\nabla_s^\perp)^2\curv_j\|_{L^2(\dd s_j)}\Big)<\infty.
    \end{equation}
    Next, choose positive diffeomorphisms $\varphi_j\colon\S^1\to\S^1$ with $|\partial_x(\gamma_j\circ\varphi_j)|\equiv v_j=\frac{\mathcal L_{\R^2}(\gamma_j)}{2\pi}$ on $\S^1$ and $\varphi_j(0)=0$ for all $j\in\N$, and denote these reparametrizations again by $\gamma_j$. Arguing as in \cite[Remark~4.1 and proof of Theorem~1.5]{schlierf2024cvpde}, these reparametrizations preserve \eqref{eq:symmetry-turningnumber}. Moreover, using \eqref{eq:euclidean-bounds}, we may up to subsequence arrange that $|\partial_x\gamma_j|\equiv v_j\to v>0$ converges in $\R$. By \eqref{eq:symmetry-turningnumber}, we have $\gamma_j(0),\gamma_j(\pi)\in \{0\}\times(0,\infty)$. Thus, by \eqref{eq:euclidean-bounds} and using $\E(\gamma_j)\leq M$, \cite[Lemma~3.7]{schlierf2025acv} yields that there exist $Z\subset \S^1$ finite, $\gamma_\infty\in W^{1,\infty}(\S^1,\R\times[0,\infty))\cap W^{2,2}_{\mathrm{loc}}(\S^1\setminus Z,\H^2)$ such that $Z=\{\gamma_\infty^{(2)}=0\}$ and
    \begin{equation}\label{eq:uniform-profile-limit}
        \gamma_j\to\gamma_\infty\text{ in $C^0(\S^1)$},\quad\gamma_j\rightharpoonup^*\gamma_\infty\text{ in $W^{1,\infty}(\S^1)$}\quad\text{and}\quad \gamma_j\rightharpoonup\gamma_\infty\text{ in $W^{2,2}_{\mathrm{loc}}(U)$}
    \end{equation}
    with $U\vcentcolon=\S^1\setminus Z=\{\gamma_\infty^{(2)}>0\}$. In particular, $|\partial_x\gamma_\infty|=v>0$ on $U$. Moreover, we get that $\gamma_\infty$ satisfies the reflection-symmetry-identity in \eqref{eq:symmetry-turningnumber} on $\S^1$ by the uniform convergence $\gamma_j\to\gamma_\infty$. Note that, in order to conclude \eqref{eq:height-limits-at-0-pi} up to reparametrizations by $x\mapsto \pm x + \pi$ or $x\mapsto \pm x$, it is sufficient to show that $\# Z=1$ and $Z\subset\{0,\pi\}$.
    
    On compact subsets of $U$, \eqref{eq:uniform-profile-limit} implies a positive lower bound for $\gamma_j^{(2)}$. By \eqref{eq:curvature-bounds} and the argument in \cite[proof of Lemma~3.4]{schlierf2024cvpde}, $(\gamma_j)_{j}$ is bounded in $W^{4,2}_{\mathrm{loc}}(U)$. Thus, we have
    \begin{equation}\label{eq:local-profile-convergence}
        \gamma_j\rightharpoonup\gamma_\infty\text{ in $W^{4,2}_{\mathrm{loc}}(U)$}\quad\text{and}\quad\gamma_j\to\gamma_\infty\text{ in $C^{3,\beta}_{\mathrm{loc}}(U)$ for every $0<\beta<1/2$}.
    \end{equation}
    So using \eqref{eq:elastic-gradient-small}, we get $\nabla\E(\gamma_\infty)=0$ a.e.\ on $U$, and \eqref{eq:geometry}, \eqref{eq:frenet} and $|\partial_x\gamma_\infty|\equiv v>0$ on $U$ imply that $\gamma_\infty$ is smooth on $U$.
    
    For the sake of contradiction, suppose that after passing to a subsequence, we have $\sup_j\mathcal L(\gamma_j)<\infty$. The estimate
    \begin{equation}\label{eq:log-height-bound}
        |\partial_s\log(\gamma^{(2)})|\leq1
    \end{equation}
    and \eqref{eq:euclidean-bounds} imply $\inf_j\min_{\S^1}\gamma_j^{(2)}>0$. Thus $U=\S^1$ and \eqref{eq:local-profile-convergence} yields $\mathcal T(\gamma_\infty)=0$ which contradicts \cite[Corollary~5.8]{muellerspener2020}. Thus, $\mathcal L(\gamma_j)\to\infty$.
    
    By \eqref{eq:euclidean-bounds} and $\mathcal L(\gamma_j)\to\infty$, we have $\#Z\geq 1$. The argument in \cite[Lemma~2.3 and proof of Theorem~1.1]{schlierf2024cvpde}, also used in \cite[Equation~(5.8)]{dallacquamuellerruppschlierf2025}, yields $\#Z\leq 2$ and
    \begin{equation}\label{eq:outer-energy-loss}
        \E(\gamma_\infty|_{\S^1\setminus Z})\leq\liminf_{j\to\infty}\E(\gamma_j)-8\#Z\leq24-\delta-8\#Z.
    \end{equation}
    For every connected component $I\subset U$, $\gamma_\infty|_I$ is complete by \eqref{eq:log-height-bound}. 
    Using \eqref{eq:symmetry-turningnumber} and \eqref{eq:outer-energy-loss}, the argument in \cite[proof of Proposition~5.4, Claim~3]{dallacquamuellerruppschlierf2025},  excludes the case $\#Z=2$ and shows that either $Z=\{0\}$ or $Z=\{\pi\}$. After reparametrizing each $\gamma_j$ by $x\mapsto \pm x+\pi$ or $x\mapsto \pm x$ if necessary, we can thus arrange
    \begin{equation}\label{eq:choice-of-Z-and-tangent-at-pi}
        Z=\{0\},\quad\gamma_\infty(0)=(0,0), \quad \partial_s^{\R^2}\!\gamma_\infty(\pi)=(-1,0) \quad\text{and}\quad\rho_\infty\vcentcolon=\gamma_\infty^{(2)}(\pi)>0.
    \end{equation}
    This proves \eqref{eq:height-limits-at-0-pi}. Since $\gamma_\infty|_{\S^1\setminus\{0\}}$ is complete and both ends converge to the origin, \Cref{lem:finite-energy-elastica,rem:complete-elastica-ends} yield that $f_{\gamma_\infty}$ is an inverted catenoid. 
    Using \eqref{eq:symmetry-turningnumber}, \eqref{eq:choice-of-Z-and-tangent-at-pi}, and \eqref{eq:local-profile-convergence}, we have 
    \begin{equation}\label{eq:left-pointing-at-pi}
        \partial_s^{\R^2}\!\gamma_j(\pi) = (-1,0)\quad\text{for all $j$ sufficiently large}
    \end{equation}
    so that \eqref{eq:horizontal-noether} implies
    \begin{equation}\label{eq:fixed-point-noether}
        P_{\gamma_j}(\pi)=-\frac{\kappa_{\gamma_j}(\pi)(\kappa_{\gamma_j}(\pi)+2)}{\gamma_j^{(2)}(\pi)}.
    \end{equation}
    Equation~\eqref{eq:local-profile-convergence} yields $\kappa_{\gamma_j}(\pi)\to\kappa_{\gamma_\infty}(\pi)=2$, so that \eqref{eq:choice-of-Z-and-tangent-at-pi} and \eqref{eq:fixed-point-noether} imply $P_{\gamma_j}(\pi) \to -\frac 8{\rho_\infty}$. 

    Finally, to obtain the local convergence of $\tilde\gamma_{j,\pi}$, note that by \eqref{eq:choice-of-Z-and-tangent-at-pi} we find that, for $s$ in a compact subset of $\R$, the image $\tilde\gamma_{j,\pi}(s)$  remains in a compact subset of $\H^2$. Since the Euclidean and hyperbolic metrics are comparable on compact subsets, the local convergence of $\tilde\gamma_{j,\pi}$ follows from \eqref{eq:local-profile-convergence}.
\end{proof}

We can now apply the quantization result \Cref{prop:kappaj-palais-smale-rigidity} to the sequence of curvatures $\tilde\kappa_{j,\pi}$ to get the following.

\begin{lemma}[Asymptotics of the hyperbolic curvatures]\label{lem:application-of-prop:kappaj-palais-smale-rigidity}
    In the setting of \Cref{lem:compactness-inverted-catenoid}, 
    we have
    \begin{equation}
        \E(\gamma_j)\to 16,\quad\kappa_{\gamma_j}(\pi)\to 2,\quad\text{and}\quad \kappa_{\gamma_j}(0)\to \pm 2,
    \end{equation}
    and, for all $(x_j)_j\subset (0,\pi)$ such that $\mathcal L(\gamma_j|_{[0,x_j]})\to\infty$ and $\mathcal L(\gamma_j|_{[x_j,\pi]})\to\infty$, we have
    \begin{equation}\label{eq:kappa-j-of-x-j-goes-to-zero}
        \kappa_{\gamma_j}(x_j)\to 0.
    \end{equation}
\end{lemma}
\begin{proof}
    \Cref{lem:compactness-inverted-catenoid} and \eqref{eq:elastic-gradient-small} yield that $\tilde\kappa_{j,\pi}$ are $L_j$-periodic with $L_j=\mathcal L(\gamma_j)\to\infty$, satisfying up to subsequence \eqref{eq:kappaj-palais-smale-assumptions} with $E\in [16,24)$, using \eqref{eq:outer-energy-loss} with $\#Z\geq 1$ for the lower bound. Thus \Cref{prop:kappaj-palais-smale-rigidity} applies with $m=2$. Moreover, by \eqref{eq:kappa-infty-pi}, we have $\tilde\kappa_{j,\pi}(0)\to 2$ so that we may choose $a_j^{(1)}=0$ for all $j\in\N$. 
    
    Since \eqref{eq:symmetry-turningnumber} yields that $\tilde\kappa_{j,\pi}$ is even for all $j\in\N$, \eqref{eq:kappaj-convergence-at-a_ji} yields for $\sigma_2\in\{-1,1\}$ that
    \begin{equation}
        \tilde\kappa_{j,\pi}(-a_j^{(2)}-s)=\tilde\kappa_{j,\pi}(a_j^{(2)}+s)\to \sigma_2\frac{2}{\cosh(s)}\quad\text{in $C^{1,\beta}_{\mathrm{loc}}(\R)$ for all $0<\beta<1/2$.}
    \end{equation}
    So $\tilde\kappa_{j,\pi}(-a_j^{(2)})\to \pm 2$, and therefore by \eqref{eq:kappaj-univ-convergence-to-zero}, \eqref{eq:aij-uniqueness} and $\mathrm{dist}_j(a_j^{(1)},-a_j^{(2)})=\mathrm{dist}_j(a_j^{(1)},a_j^{(2)})\to\infty$, we obtain 
    \begin{equation}
        2\operatorname{dist}_j(\mathcal L(\gamma_j)/2,a_j^{(2)})\leq \operatorname{dist}_j(a_j^{(2)},-a_j^{(2)})\to 0,
    \end{equation} 
    which shows that we may take $a_j^{(2)}=\mathcal L(\gamma_j)/2$. By \eqref{eq:kappaj-convergence-at-a_ji}, we thus have that
    \begin{equation}
        \kappa_{\gamma_j}(0)=\tilde\kappa_{j,\pi}(\mathcal L(\gamma_j)/2) \to \sigma_2\cdot 2.
    \end{equation}
    Moreover, \eqref{eq:kappa-j-of-x-j-goes-to-zero} follows from \eqref{eq:kappaj-univ-convergence-to-zero}. 
\end{proof}

The main new ingredient is the following lemma which focuses on the asymptotics of the mean curvature $H$ at the points where the height $\gamma_j^{(2)}\to 0$ vanishes: at these points, we do not have an analog of \eqref{eq:gamma-infty-pi}. Instead, the new idea is to prove that the Noether quantity $P_{\gamma_j}$ from \eqref{eq:horizontal-noether}, whose limit at $x=\pi$ is known by \eqref{eq:kappa-infty-pi}, actually converges \emph{uniformly} on $\S^1$ to the same constant. This uniform limit then gives the missing link between the known outer limit from \Cref{lem:compactness-inverted-catenoid} and the unknown behavior close to the concentration point $x=0$. 

The classical approach would be to \emph{classify the blow-up} at $x=0$, i.e., the limit of $\gamma_j/\gamma_j^{(2)}(0)$. This limit can be shown to be a catenoid having zero mean curvature. However, our finer analysis gives additional information about the limiting mean curvature before blow-up, see \eqref{eq:hinfty-and-kappainfty-at-0} below. Due to the scaling behavior of the mean curvature, this information is necessarily lost in the blow-up limit.

\begin{lemma}\label{lem:Astranomical_contribution}
    In the setting of both \Cref{lem:compactness-inverted-catenoid,lem:application-of-prop:kappaj-palais-smale-rigidity}, we have
    \begin{align}
        P_{\gamma_j}\to -\frac{8}{\rho_\infty}\quad\text{in $C^0(\S^1)$}. \label{eq:unif-conv-Pgammaj}
    \end{align}
    Fix the subsequence, reparametrizations and choice of $\tilde\gamma_{j,\pi}$ as in \Cref{lem:compactness-inverted-catenoid}. Then define
    \begin{equation}\label{eq:def-tildegamma-j-0}
        \tilde\gamma_{j}(s)=\tilde\gamma_{j,\pi}(s-\frac12\mathcal L(\gamma_j))\quad\text{for all $s\in\R$}.
    \end{equation}
    That is, $\tilde\gamma_{j}$ is the $\mathcal L(\gamma_j)$-periodic lift to $\R$ of the reparametrization of $\gamma_j$ by hyperbolic arclength fixing the value at $s=x=0$, and thus satisfies \eqref{eq:symmetry-turningnumber} for each $j\in\N$. Moreover, for the curvatures $\tilde\kappa_{j}$ of $\tilde\gamma_{j}$ and mean curvatures $\tilde H_{j}$ of $f_{\tilde\gamma_{j}}$, we have
    \begin{equation}\label{eq:hinfty-and-kappainfty-at-0}
        \tilde\kappa_{j} \to -\frac{2}{\cosh(\cdot)} \quad\text{and}\quad \tilde H_{j}\to -\tilde H_\infty \quad \text{in $C^{1,\beta}_{\mathrm{loc}}(\R)$ for every $0<\beta<1/2$}
    \end{equation}
    with $\tilde H_\infty$ as in \eqref{eq:partc-Hinfty}.
\end{lemma}
\begin{proof}
    Using \eqref{eq:horizontal-noether_estimate}, the fundamental theorem of calculus, Cauchy--Schwarz, \eqref{eq:gradient-norm}, and \eqref{eq:area-measure}, we obtain
    \begin{align}
        \max_{\S^1}P_{\gamma_j}-\min_{\S^1}P_{\gamma_j} &\leq\int_{\S^1}\frac{|\nabla\E(\gamma_j)|_g}{\gamma_j^{(2)}}\dd s_j \leq \Big(\int_{\S^1}\frac{|\nabla\E(\gamma_j)|_g^2}{(\gamma_j^{(2)})^4}\dd s_j\Big)^{1/2}\Big(\int_{\S^1}(\gamma_j^{(2)})^2\dd s_j\Big)^{1/2} \\
        &=\frac{2}{\pi}\|\nabla\W(f_{\gamma_j})\|_{L^2(\dd\mu_{f_{\gamma_j}})}\to0\label{eq:horizontal-noether-oscillation}
    \end{align}
    by \eqref{eq:compactness-assumptions}. Together with the pointwise convergence $P_{\gamma_j}(\pi)\to -\frac{8}{\rho_\infty}$ from \eqref{eq:kappa-infty-pi}, we get \eqref{eq:unif-conv-Pgammaj}.
    
    Next, we note that \eqref{eq:symmetry-turningnumber} and \eqref{eq:left-pointing-at-pi} yield for all sufficiently large $j$ that
    \begin{equation}\label{eq:choice-orientation-0-pi}
        \partial_s^{\R^2}\!\gamma_j(0)=\partial_s^{\R^2}\!\gamma_j(\pi)=(-1,0).
    \end{equation}
    Combined with \eqref{eq:height-limits-at-0-pi}, \eqref{eq:horizontal-noether}, and \eqref{eq:unif-conv-Pgammaj}, we get 
    \begin{equation}\label{eq:application-Pgammaj-convergence}
        \kappa_{\gamma_j}(0)(\kappa_{\gamma_j}(0)+2) = -\gamma_j^{(2)}(0) P_{\gamma_j}(0) \to 0.
    \end{equation}
    In particular, using $\tilde\kappa_{j}(0)=\kappa_{\gamma_j}(0)$, we get that $\lim_{j\to\infty}\tilde\kappa_{j}(0)=\lim_{j\to\infty}\kappa_{\gamma_j}(0)\in \{0,-2\}$. By \Cref{lem:application-of-prop:kappaj-palais-smale-rigidity}, this implies $\lim_{j\to\infty}\tilde\kappa_{j}(0)=-2$ and $\tilde\kappa_{j}\to -2/\cosh$ in $C^{1,\beta}_{\mathrm{loc}}(\R)$, using \eqref{eq:kappaj-convergence-at-a_ji} and \eqref{eq:def-tildegamma-j-0}.
     
    Using \eqref{eq:symmetry-turningnumber} and \eqref{eq:H-A0} we have that $\tilde H_j$ is even, thus $\partial_s\tilde H_j(0)=0$. Combined with \eqref{eq:horizontal-noether-Hformula} and \eqref{eq:choice-orientation-0-pi}, we get that
    \begin{equation}
        P_{\gamma_j}(0)=P_{\tilde\gamma_{j}}(0)=-4\tilde H_j(0)(\tilde \gamma_{j}^{(2)}(0)\tilde H_j(0)+1).
    \end{equation}
    Since $\tilde \gamma_{j}^{(2)}(0)\tilde H_j(0)+1=-\tilde\kappa_{j}(0)/2=-\kappa_{\gamma_j}(0)/2$ by \eqref{eq:H-A0} and \eqref{eq:choice-orientation-0-pi}, we get
    \begin{equation}\label{eq:mean-curvature-initial-data}
        \tilde H_j(0)=\frac{P_{\gamma_j}(0)}{2\kappa_{\gamma_j}(0)}\to\frac{2}{\rho_\infty}\quad\text{and}\quad\partial_s\tilde H_j(0)=0.
    \end{equation}
    Using $\mathcal L(\gamma_j)\to\infty$, \eqref{eq:forced-Jacobi-equation-H}, Cauchy--Schwarz, \eqref{eq:gradient-norm}, and \eqref{eq:euclidean-bounds}, we obtain for any compact interval $I\subset\R$ and $j$ sufficiently large that
    \begin{align}
        \int_{I}|2\partial_s^2\tilde H_j+\tilde\kappa_{j}^2\tilde H_j|^2\dd s &\leq
        \int_{\S^1}|2\partial_s^2H_{f_{\gamma_j}}+\kappa_{\gamma_j}^2H_{f_{\gamma_j}}|^2\dd s_j \\
        &\leq\frac14\int_{\S^1}\frac{|\nabla\E(\gamma_j)|_g^2}{(\gamma_j^{(2)})^2}\dd s_j \leq C\|\nabla\W(f_{\gamma_j})\|_{L^2(\dd\mu_{f_{\gamma_j}})}^2\to0.
    \end{align}
    Together with \eqref{eq:mean-curvature-initial-data} and the uniform convergence of $\tilde\kappa_{j}$, Gronwall's inequality applied to \eqref{eq:forced-Jacobi-equation-H} (written as a first order system) yields local $W^{2,2}$-bounds for $\tilde H_j$. By compactness and passing to the limit in \eqref{eq:forced-Jacobi-equation-H}, every subsequential limit $\hat H_{\infty}\in W^{2,2}_{\mathrm{loc}}$ satisfies 
    \begin{equation}
        2\partial_s^2\hat H_\infty+\frac{4}{\cosh^2(s)}\hat H_\infty=0\text{ in $L^2_{\mathrm{loc}}$},\quad\hat H_\infty(0)=\frac{2}{\rho_\infty}\quad\text{and}\quad\partial_s\hat H_\infty(0)=0.
    \end{equation}
    The unique solution is $\hat H_{\infty}(s)=\frac{2}{\rho_\infty}(1-s\tanh(s))$, which implies \eqref{eq:hinfty-and-kappainfty-at-0}.
\end{proof}

\section{Proof of the main results}

Our first main result establishes \eqref{eq:intro-velocity_estimate}. 

\begin{theorem}\label{thm:willgrad-lower-bound}
    For every $\delta\in(0,4\pi)$, there exists $c(\delta)>0$ with the following property. Suppose $f_\gamma$ is a torus of revolution such that its profile curve $\gamma$ satisfies \eqref{eq:symmetry-turningnumber}.
    If $\W(f_\gamma)\leq 12\pi-\delta$, then
    \begin{equation}
        \A(f_\gamma)^2 \int_{\T^2} |\nabla\W(f_\gamma)|^2\dd\mu_{f_\gamma} \geq c(\delta).
    \end{equation}
\end{theorem}
\begin{proof}
    For the sake of contradiction, suppose that the claim fails for some $\delta\in(0,4\pi)$. Since $\W(f_\gamma)$ and $\A(f_\gamma)^2\|\nabla\W(f_\gamma)\|_{L^2(\dd\mu_{f_\gamma})}^2$ are invariant under dilations, using \eqref{eq:BG} there exists a sequence of smooth immersions $\gamma_j\colon\S^1\to\H^2$ satisfying \eqref{eq:symmetry-turningnumber} and \eqref{eq:compactness-assumptions}. We apply \Cref{lem:compactness-inverted-catenoid,lem:application-of-prop:kappaj-palais-smale-rigidity,lem:Astranomical_contribution} and, after passing to a subsequence, use their notation and the parametrizations, shift, and orientation.
    Fix $R>0$ such that $R\tanh(R)>1$. 
    By \eqref{eq:hinfty-and-kappainfty-at-0}, \eqref{eq:def-tildegamma-j-0} and \eqref{eq:partc-Hinfty}, we obtain
    \begin{equation}
        \tilde H_j(R) \to-\tilde H_\infty(R)<0\quad\text{and}\quad\tilde H_j(\mathcal L(\gamma_j)/2-R)=\tilde H_{j,\pi}(-R) \to\tilde H_\infty(-R)=\tilde H_\infty(R)>0.\label{eq:willgrad-H-signs}
    \end{equation}
    Moreover, \eqref{eq:unif-conv-Pgammaj} implies
    \begin{equation}\label{eq:unif-conv-P-tildegamma-j-0}
        \Big\|P_{\tilde\gamma_{j}}+\frac{8}{\rho_\infty}\Big\|_{L^\infty(\R)}=\Big\|P_{\gamma_j}+\frac{8}{\rho_\infty}\Big\|_{L^\infty(\S^1)}\to0.
    \end{equation}
    Since $\mathcal L(\gamma_j)\to\infty$, for all sufficiently large $j$ we have $\mathcal L(\gamma_j)>4R$ and
    \begin{equation}
        \tilde H_j(R)<0<\tilde H_j(\mathcal L(\gamma_j)/2-R)\quad\text{and}\quad P_{\tilde\gamma_{j}}<0\quad\text{on $[R,\mathcal L(\gamma_j)/2-R]$}.
    \end{equation}
    In particular, there exists $x_j'\in (R,\mathcal L(\gamma_j)/2-R)$ such that $\tilde H_j(x_j')=0$ and $\partial_s\tilde H_j(x_j')\geq 0$. Using \eqref{eq:hinfty-and-kappainfty-at-0}, it follows that 
    \begin{equation}\label{eq:xj-prime-goes-to-infty}
        \mathcal L(\tilde\gamma_{j}|_{[0,x_j']}) = x_j'\to\infty.
    \end{equation}
    Evaluating \eqref{eq:horizontal-noether-Hformula} at $x_j'$ and using \eqref{eq:unif-conv-P-tildegamma-j-0} yields that
    \begin{equation}\label{eq:negative-noether-quantity}
        0 > \frac{P_{\tilde\gamma_{j}}(x_j')}{4} = \partial_s\tilde H_j(x_j') \cdot \partial_{s}^{\R^2}\tilde\gamma_{j}^{(2)}(x_j'),
    \end{equation}
    that is, $\partial_x\tilde\gamma_{j}^{(2)}(x_j')<0$. Using \eqref{eq:height-limits-at-0-pi}, it follows that, for all $j$ sufficiently large, there exists $x_j\in (0,\pi)$ a local minimum of $\gamma_j^{(2)}$ with $\mathcal{L}(\gamma_j\vert_{[0,x_j]})\geq x_j'\to\infty$. From \eqref{eq:gamma-infty-pi}, we can see that $\mathcal L(\gamma_j|_{[x_j,\pi]})\to\infty$. Since $x_j$ is a local minimum of $\gamma_j^{(2)}$, we have $\partial_s\gamma_j^{(2)}(x_j)=0$, $\partial_s^2\gamma_j^{(2)}(x_j)\geq 0$ and $(\partial_s\gamma_j^{(1)}(x_j))^2=|\partial_s\gamma_j(x_j)|^2=(\gamma_j^{(2)}(x_j))^2$. Thus, \eqref{eq:cov-der-H2} implies
    \begin{equation}
        (\nabla_s\partial_s\gamma_j)^{(2)}(x_j) = \partial_s^2\gamma_j^{(2)}(x_j)+\frac{\partial_s\gamma_j^{(1)}(x_j)}{\gamma_j^{(2)}(x_j)}\vec n_{\gamma_j}^{(2)}(x_j) = \partial_s^2\gamma_j^{(2)}(x_j)+\gamma_j^{(2)}(x_j)\geq \gamma_j^{(2)}(x_j),
    \end{equation}
    so that \eqref{eq:frenet} yields $|\kappa_{\gamma_j}(x_j)|=|\nabla_s\partial_s\gamma_j(x_j)|_g\geq 1$, contradicting \eqref{eq:kappa-j-of-x-j-goes-to-zero}.
\end{proof}

Note that the conformal class degenerates, and the conformal factor necessarily needs to blow up along the degenerating Palais--Smale sequences we consider. In particular, the results in \cite{bernardrivere2011} do not apply. 

\begin{remark}[Spheres of revolution]\label{rem:spheres_not_good}
    The argument used in the proof of \Cref{thm:willgrad-lower-bound} fails when working with Palais--Smale sequences of spheres of revolution such as the ones considered in \cite[Theorem~5.1, Figure~5.1]{blatt2009}. In this case, the macroscopic \emph{outer limit} $\gamma_\infty$ would consist of two round spheres touching at one point, contrary to the inverted catenoid we obtain in \Cref{lem:compactness-inverted-catenoid}. For round spheres, the Noether quantity $P_{\gamma_\infty}$ vanishes. Therefore, the argument in \eqref{eq:negative-noether-quantity} showing the existence of further minima of $\gamma^{(2)}$ and thus leading to a contradiction fails. A similar problem occurs for figure-eight tori of revolution when dropping the reflection symmetry assumption \eqref{eq:symmetry-turningnumber}. Indeed, also in this case, the macroscopic limit can consist of two round spheres.
\end{remark}

We next prove a quantitative version of \Cref{thm:intro-main}.

\begin{theorem}\label{thm:main}
    For every $\delta\in(0,4\pi)$, there exists $0<T(\delta)<\infty$ with the following property. Suppose $f_0$ is a torus of revolution with
    \begin{equation}
        \W(f_0)\leq 12\pi-\delta
    \end{equation}
    such that its profile curve $\gamma_0$ satisfies \eqref{eq:symmetry-turningnumber}. 
    Then the maximal Willmore flow $f\colon[0,T_{\max})\times\T^2\to\R^3$ with initial datum $f_0$ satisfies
    \begin{equation}
        T_{\max}\leq \A(f_0)^2 T(\delta).
    \end{equation}
\end{theorem}
\begin{proof}
    Suppose that $f_0=f_{\gamma_0}$ is a torus of revolution as in the statement. Using \cite[Lemma~3.3]{dallacquamullerschatzlespener2020}, we have that $f(t)=f_{\gamma(t)}$ for a smooth family of immersions $\gamma\colon[0,T_{\max})\times\S^1\to\H^2$. Moreover, we have $\W(f(t))\leq \W(f_0)\leq 12\pi-\delta$ for all $t\in[0,T_{\max})$ and, arguing as in \cite[proof of Proposition~5.4]{dallacquamuellerruppschlierf2025}, $\gamma(t)$ satisfies \eqref{eq:symmetry-turningnumber} for all $t\in[0,T_{\max})$. Moreover, \cite[proof of Lemma~D.6]{dallacquamullerschatzlespener2020} yields that
    \begin{equation}
        \A(f(t))\leq \A(f_0)+2\W(f_0) t^\frac{1}{2}\quad\text{for all $0\leq t<T_{\max}$}.
    \end{equation}
    Fixing the constant $c(\delta)$ from \Cref{thm:willgrad-lower-bound}, we get
    \begin{align}
        -\frac{\dd}{\dd t}\W(f(t)) = \int_{\T^2} |\nabla\W(f(t))|^2\dd\mu_{f(t)} \geq c(\delta) \A(f(t))^{-2} \geq c(\delta) (\A(f_0)+2\W(f_0)t^{\frac12})^{-2}. 
    \end{align}
    Using $\W(f_0)\leq12\pi-\delta$ and the fact that $\W(f(t)) \geq 8\pi$ by the Li--Yau inequality and $\mathcal T(\gamma(t))=0$ for all $t\in[0,T_{\max})$, integrating on $[0,T_{\max})$ yields 
    \begin{equation}
        4\pi-\delta\geq \W(f_0)-\lim_{t\to T_{\max}}\W(f(t))\geq \frac{c(\delta)}{2\W(f_0)^2} F(2\W(f_0)/\A(f_0)\cdot \sqrt{T_{\max}})
    \end{equation}
    where $F(x)=\log(1+x)+\frac{1}{1+x}-1$. This yields the claim with $T(\delta)$ depending only on $c(\delta)$ and $\delta$.
\end{proof}

Finally, we prove \Cref{thm:intro-main}.

\begin{proof}[Proof of \Cref{thm:intro-main}]
    Suppose that $\W(f_0)=12\pi$ and $\partial_{t,0}\W(f(t))=0$. Writing $f_0=f_{\gamma_0}$,  $\gamma_0$ is an elastica with turning number zero by \Cref{prop:elastic-willmore-identities}, contradicting \cite[Corollary~5.8]{muellerspener2020}. Thus $\partial_{t,0}\W(f(t))<0$ and \Cref{thm:main} applies with $f_0$ replaced by $f(t)$ for $t>0$ sufficiently small.
\end{proof}

\appendix
\crefalias{section}{appendix}%
\crefalias{subsection}{appendix}

\section*{AI declaration}
Some of the ideas in this work were developed interactively with the help of OpenAI's GPT-6 Astra Ultra model, which was given an incomplete and unpublished draft by the authors on the topic on November 2, 2026. Its main contribution is the argument in \Cref{lem:Astranomical_contribution}, which completed our earlier line of reasoning. The authors have checked, polished, and simplified this argument, and incorporated it into their previous ideas. The manuscript was written independently by the authors, who take full responsibility for its content. Claude Opus 5.5 by Anthropic was used to create the TikZ code for \Cref{fig:introduction} from hand-drawn notes and to edit the wording of this declaration.
\printbibliography
\end{document}